\documentstyle[12pt,amssymb,amsmath,epsf]{article}

\def\R {{\Bbb R}}

\def \Z {{\Bbb Z}}
\def \S {{\Bbb S}}
\def\rth{{\Bbb R}^3}

\newtheorem{theorem}{Theorem}
\newtheorem{proposition}{Proposition}
\newtheorem{definition}{Definition}
\newtheorem{lemma}{Lemma}

\newtheorem{remark}{Remark}
\newtheorem{claim}{Claim}

\newcommand{\myskip}[1]{}

\begin{document}
\title{Saddle towers with infinitely many ends}
\author{Laurent Mazet\\M. Magdalena Rodr\'\i guez\thanks{Research partially supported by grants from R\'egion Ile-de-France and a MEC/FEDER grant no. MTM2004-02746.}\\Martin Traizet}
\date{April 24, 2006}

\maketitle

\section{Introduction}\label{secintro}
The topic of minimal surfaces in flat 3-manifolds, with finite genus but infinite
total curvature, has recently attracted some attention~\cite{hauswirth-pacard,howe3}. 
In the complete flat 3-manifold $\R^2\times\S^1$, the only known examples of properly embedded minimal surfaces with infinite total curvature come from doubly or triply periodic minimal surfaces in $\R^3$. 
In particular, they are all periodic in $\R^2\times \S^1$. 
In this paper, we point out an application of a theorem by Jenkins and Serrin~\cite{jes1}
to construct properly embedded minimal surfaces in $\R^2\times\S^1$ with genus zero and infinite total
curvature. We~prove:

\begin{theorem}\label{mainth}
There exists a properly embedded singly periodic minimal surface ${\cal M}$ in $\R^3$ with bounded (Gaussian) curvature, whose quotient by all its periods has genus zero, infinitely many ends and exactly one limit end.
\end{theorem}

Recall that Scherk's singly periodic minimal surface can be constructed as follows:
consider the unit square, and mark its two horizontal edges with $+\infty$ and its two vertical edges with $-\infty$.
By the theorem of Jenkins and Serrin~\cite{jes1}, there exists a function $u$ which solves the Jenkins-Serrin problem on the
square, namely, whose graph is minimal in the interior of the square, and which goes to $\pm \infty$ on the edges, as indicated by the marking.
The graph of $u$ is bounded by four vertical lines above the vertices of the square and is a fundamental piece for Scherk's doubly periodic minimal surface.
The conjugate minimal surface 
of the graph of $u$ is bounded by four horizontal symmetry curves, lying in two horizontal planes at distance 1 from each other.
By reflecting about one of the two symmetry planes, we obtain a fundamental domain for Scherk's singly periodic minimal surface,
which has period $T=(0,0,2)$, and four ends in the quotient.

\medskip

H. Karcher~\cite{ka4} has generalized this construction by replacing the unit square by any convex polygonal domain $\Omega$
with $2k$ edges of length one, $k\geq 2$.
To satisfy the hypothesis of the theorem of Jenkins and Serrin, the domain $\Omega$ must be assumed to be non-special,
see definition~\ref{def1} below 
(this known fact does not seem to have been written yet, so we provide a proof in the appendix).
Solving the Jenkins-Serrin problem on $\Omega$, taking the conjugate and reflecting,
one obtains a properly embedded singly periodic minimal surface with period $T=(0,0,2)$,
$2k$ Scherk-type ends and genus zero in the quotient.
These surfaces are now called Karcher's Saddle Towers, 
and have recently been classified as the only properly embedded singly
periodic minimal surfaces in $\R^3$ with genus zero and finitely many Scherk-type ends in the quotient~\cite{PeTra1}.
\begin{definition}
\label{def1}
We say a convex polygonal domain with $2k$ unitary edges is special if $k\geq 3$ and its boundary is a parallelogram with two sides of length one and two sides of length $k-1$.
\end{definition}

\medskip

In this paper we follow the same strategy except that we start with an unbounded convex domain $\Omega$ with infinitely many edges,
so we end up with a minimal surface with infinitely many ends, as desired.
More precisely, we consider an unbounded convex domain $\Omega\subset\R^2$ such that:
\begin{enumerate}
\item The boundary $\partial\Omega$ of $\Omega$ is a polygonal curve with an infinite number of edges, all of
length one.
\item $\Omega$ is neither the plane, nor a half plane, nor a strip, nor an infinite special domain, see definition~\ref{def2} below. 
\end{enumerate}

\begin{definition}
\label{def2}
An unbounded convex polygonal domain is said to be special when its boundary is made of two parallel half lines and one edge of length one (such a domain may be seen as a limit of special domains with $2k$ edges, when
$k\to\infty$). 
\end{definition}

Given a domain $\Omega$ as above, mark the edges on its boundary alternately by $+\infty$ and $-\infty$. 
In section~\ref{secJSproblem}, we solve the Jenkins-Serrin problem for $\Omega$. 
In order to do this, we consider an exhaustion of $\Omega$ by bounded convex domains $\Omega_n$ and solve the Jenkins-Serrin problem on each $\Omega_n$, obtaining a solution $u_n$ in $\Omega_n$. 
Then we prove that the sequence $\{u_n\}_n$ has a limit $u$.
Such a function $u$, which is defined on $\Omega$, has the required behavior on the boundary and its
graph $M$ is minimal. Taking the conjugate minimal surface of $M$ and extending by symmetry, we obtain the desired
minimal surface. Such a surface can be seen as a limit, when $k\to\infty$, of a sequence of Karcher's Saddle Towers
with $2k$ ends.
In section~\ref{secBounded}, we prove this surface has bounded curvature.
Finally, in section~\ref{secAsymptotic} we study the asymptotic behavior of this surface.

\section{Preliminaries}\label{secprelim}
Let $u=u(x_1,x_2)$ be a solution of the minimal graph equation:
\begin{equation}\label{eqmin}
(1+u_{2}^2) u_{11}-2u_{1} u_{2} u_{12}+(1+u_{1}^2) u_{22}=0,
\end{equation}
defined on a simply-connected domain $D\subset\R^2$. 
By an elementary computation,
\begin{equation}
\label{conjugate}
d\psi_u:=\frac{u_{1}}{\sqrt{1+|\nabla u|^2}}\, d x_2-\frac{u_{2}}{\sqrt{1+|\nabla u|^2}}\, d x_1
\end{equation}
is an exact form in $D$. Hence there exists a function $\psi_u=\psi_u(x_1,x_2)$, called {\it conjugate function of $u$},
whose differential is given by (\ref{conjugate}).
Note that $\psi_u$ is well defined up to an additive constant.
In fact, if we write $X(x_1,x_2)=(x_1,x_2,u(x_1,x_2))$ and call $X^*=X^*(x_1,x_2)$ its
conjugate minimal immersion, then the third coordinate function of $X^*$ can be written as $X^*_3(x_1,x_2)=\psi_u(x_1,x_2)$ (although the conjugate surface is not the graph of $\psi_u$).

Since $|\nabla\psi_u|=\frac{|\nabla u|}{\sqrt{1+|\nabla u|^2}}<1$, $\psi_u$ is a Lipschitz function, so
it can be extended continuously to $\partial D$.

Next we expose some results related to the convergence of a sequence $\{u_n\}_n$ of minimal graphs defined on $D$.
They are based on the theory developed by L. Mazet~\cite{mazet1,mazet2}, following the ideas of Jenkins and Serrin 
(the main improvement over the work of Jenkins and Serrin is that we do not require monotonicity of the sequence $\{u_n\}_n$).

Given a sequence $\{u_n\}_n$ of solutions for the minimal graph equation in~$D$, 
define the {\it convergence domain} of the sequence $\{u_n\}_n$ as
\[
{\cal B}(u_n)=\left\{p\in D\ |\ \{|\nabla u_n(p)|\}_n \mbox{ is bounded}\right\}.
\]
For each component $D'$ of ${\cal B}(u_n)$, there is a subsequence of $\{u_n-u_n(Q)\}_n$ converging uniformly
on compact sets of $D'$ to a solution of (\ref{eqmin}), where $Q$ is some fixed point of $D'$.
This fact justifies the name for ${\cal B}(u_n)$.
Moreover, 
\[
D-{\cal B}(u_n)=\cup_{i\in I} L_i ,
\]
where each $L_i\subset D$ is a component of the intersection of a straight line with~$D$, for each $i\in I$.
The straight lines $L_i$ are called {\it divergence lines}.

Clearly, to ensure the convergence of a subsequence of $\{u_n\}_n$ on $D$, it suffices to prove there are no divergence lines. The following lemmas~\ref{lem1} and~\ref{lem2} can be useful to conclude this. 

\begin{lemma}[\cite{mazet2}]
\label{lem1}
If $T\subset\partial D$ is an open straight segment such that each $u_n$ diverges to $+\infty$ when we approach
$T$, then a divergence line cannot end in~$T$.
\end{lemma}

\begin{lemma}[\cite{mazet1}]
\label{lem2}
Given a segment $T$ contained in a divergence line, it holds ${\int_T d\psi_n\to \pm |T|}$.
\end{lemma}

Once we have ensured the convergence of the sequence $\{u_n\}_n$ to a solution $u$ of the minimal graph equation, 
the next natural step is to understand the behavior of $u$ on the boundary of $D$.

\begin{lemma}[\cite{jes1,mazet2}]
\label{lem4}
Let $u$ be a solution of (\ref{eqmin}) on $D$, and $T\subset\partial D$ be an open straight segment 
oriented as $\partial D$. 
Then, $\int_T d\psi_u= |T| $ if and only if $u$ diverges to $+\infty$ on $T$.
\end{lemma}

In section \ref{secBounded} and \ref{secAsymptotic},
we will consider sequences $\{u_n\}_n$ of solutions for
the minimal graph equation defined on domains $D_n$ which
are not fixed.
For a sequence of convex domains $D_n$ in $\R^2$, we define its limit domain $D_{\infty}$ as the set of points of $\R^2$ that admit a neighborhood contained in every $D_n$, for $n$ large enough.
$D_\infty$ is a convex open set.
Given a solution $u_n$ of the minimal surface equation on $D_n$,
we defined the convergence domain ${\cal B}(u_n)$ of the sequence $\{u_n\}_n$, 
 as the set of points $p\in D_{\infty}$ such that the sequence $\{|\nabla u_n(p)|\}_n$ (which is defined for
$n$ large enough) is bounded.
As above, the complement of the convergence domain is a union of divergence
lines in $D_{\infty}$.
We need a generalization of lemmas \ref{lem1} and \ref{lem2} for this setting.

Let $\psi_n$ be the conjugate function of $u_n$.
Since $\psi_n$ is Lipschitz, a subsequence of $\{\psi_n\}_n$ converges 
to a Lipschitz function $\psi_\infty$ on $D_\infty$. 
Since $\psi_\infty$ is Lipschitz, it extends to the boundary of $D_\infty$.

\begin{lemma}
\label{laurent}
Let $\{D_n\}_n$  be the sequence of domains defined by $D_n=\{(x,y)\in(0,1)^2\,|\, y<a_n x+b_n \}$. 
We assume that $a_n\rightarrow 0$ and $b_n\rightarrow 1/2$; \textit{i.e.} $\{D_n\}_n$ converges to
$D_\infty=(0,1)\times(0,1/2)$. 
Let $u_n$ be a solution of (\ref{eqmin}) on $D_n$ such that $u_n=+\infty$ on the segment $\Gamma_n=\{y=a_nx+b_n\}\cap \partial D_n$, $\psi_n$ be the conjugate function of $u_n$, and $\psi_\infty$ be a limit of $\{\psi_n\}_n$ on $D_\infty$. 
Then:
\begin{itemize}
\item[(i)] No divergence line of $\{u_n\}_n$ ends at $\Gamma_\infty=(0,1)\times\{1/2\}$. 
\item[(ii)] $\psi_\infty(s,1/2)-\psi_\infty(t,1/2)=t-s$ for every $0\leq s<t\leq 1$.
\end{itemize}
\end{lemma}
\begin{proof}
Let $A_n=(0,b_n)$ and $B_n=(1,a_n+b_n)$ be the end-points of $\Gamma_n$. 
By hypothesis, $A_n\rightarrow A_\infty=(0,1/2)$ and $B_n\rightarrow B_\infty=(1,1/2)$. 
Since $u_n$ takes the value $+\infty$ along $\Gamma_n$, Lemma \ref{lem4} says that
$\psi_n(A_n)-\psi_n(B_n)=|A_nB_n|$. 
By other hand, each $\psi_n$ is Lipschitz, so $\psi_n(A_n)-\psi_n(B_n)\rightarrow \psi_\infty(A_\infty)-\psi_\infty(B_\infty)$. 
Thus $\psi_\infty(A_\infty)-\psi_\infty(B_\infty)=|A_\infty B_\infty|=1$.
Since $\psi_\infty$ is Lipschitz, we conclude item {\it (ii)}. 

Suppose that $\{u_n\}_n$ has a divergence line $L$ which ends at a point $P\in\Gamma_\infty$, 
and  consider two points $P',Q$ in $L$ such that $P'$ is between $P$ and $Q$. 
For $n$ big enough, $P'\in D_n$, and $\psi_n(P')-\psi_n(Q)\rightarrow \psi_\infty(P')-\psi_\infty(Q)$. 
But since $L$ is a divergence line,
$|\psi_n(P')-\psi_n(Q)|\rightarrow |P'Q|$ (Lemma \ref{lem2}), thus
$|\psi_\infty(P')-\psi_\infty(Q)|=|P'Q|$. 
Letting $P'$ tend to $P$, we get $|\psi_\infty(P)-\psi_\infty(Q)|=|PQ|$. 
We can assume $\psi_\infty(P)-\psi_\infty(Q)=|PQ|$ 
(if it equals to $-|PQ|$, we consider $B_\infty$ instead of $A_\infty$ in following argument). 
Using {\it (ii)} and the fact that $\psi_\infty$ is Lipschitz in $D_\infty$, we have:
\begin{align*}
|QP|+|PA_\infty|&= \left(\psi_\infty(P)-\psi_\infty(Q)\right)+
\left(\psi_\infty(A_\infty)-\psi_\infty(P)\right)\\
&=\psi_\infty(A_\infty)-\psi_\infty(Q)\le |A_\infty Q|.
\end{align*}
This contradicts the triangle inequality and proves Lemma~\ref{laurent}. 
\end{proof}

\medskip
Finally, we have the following uniqueness result for the limit $u$ under some constraints.
\begin{lemma}[\cite{mazet3}]
\label{lem5}
Let $u_1$ and $u_2$ be two solutions of (\ref{eqmin}) in a connected domain~$D$, whose conjugate functions $\psi_{u_1},\psi_{u_2}$ are bounded in $D$ and satisfy $\psi_{u_1}=\psi_{u_2}$ on $\partial D$.
Then $u_1-u_2$ is constant in $D$.
\end{lemma}

\section{Solving the Jenkins-Serrin problem on $\Omega$}
\label{secJSproblem}
Let $\Omega$ be an unbounded convex domain as in the introduction.
We choose a vertex $p_0$ such that the inner angle at $p_0$ is less than $\pi$.
We label the vertices $p_i$, $i\in\Z$, in the order that we meet them when traveling along
the boundary of $\Omega$ with its natural orientation.
We mark the edge $(p_i,p_{i+1})$ with $+\infty$ if $i$ is even and $-\infty$ if $i$ is odd.

\begin{proposition}
\label{prop1}
The Jenkins-Serrin problem on $\Omega$ has a solution $u$ satisfying $0\leq\psi_u\leq 1$ 
($u$ is the unique solution to the Jenkins-Serrin problem with bounded
conjugate function on $\Omega$).
\end{proposition}
\begin{proof}
Given $n\geq 1$, the chord $[p_{-n},p_n]$ divides $\Omega$ into two components. 
Call $U_n$ the bounded one. Let $\Omega_n$ be the union of $U_n$ and its symmetric image about the
midpoint of the segment $[p_{-n},p_n]$. We also extend by symmetry the marking on the edges.
Since $\Omega$ is an unbounded convex domain, the sum of the inner angles of $U_n$
at $p_{-n}$ and $p_n$ is at most $\pi$. This implies that $\Omega_n$ is a (bounded) convex domain.

Let us prove that $\Omega_n$ is non-special.
If $n=1$ then this is true by definition (a special domain has at least six edges).
Assume that $n\geq 2$.
If $\Omega_n$ were special, then $p_0$ would be a corner of the parallelogram $\Omega_n$ because of
the way we chose it. Then either $p_{-2}p_{-1}p_0p_1$ or $p_{-1}p_0p_1p_2$ would be a rhombus.
Since $\Omega$ is an unbounded convex domain, it follows that $\Omega$ would be an infinite special
domain, a contradiction.
\begin{figure}\begin{center}
\epsfysize=5cm \epsffile{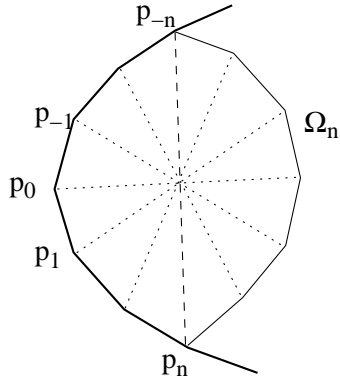}
\end{center}\caption{Definition of the domain $\Omega_n$}
\end{figure}

Hence $\Omega_n$ is non-special, so by Proposition~\ref{prop-special} that will be proven in the appendix, 
it satisfies the hypotheses of Jenkins and Serrin. 
Let $u_n$ be the solution to the Jenkins-Serrin problem on $\Omega_n$ normalized by $u_n(Q)=0$, 
where $Q$ is some fixed point in $\Omega_1$.  
Denote by $\psi_n$ the conjugate function associated to $u_n$, normalized so that $\psi_n(p_0)=0$.
From Lemma~\ref{lem4} we have $\int_{p_i}^{p_{i+1}} d\psi_n=(-1)^i$, which implies that $\psi_n(p_i)$ is equal
to $0$ if $i$ is even, and equal to $1$ if $i$ is odd. Moreover, $\psi_n$ is an affine function on each edge,
so $0\leq\psi_n\leq 1$ on $\partial\Omega_n$.
Since the domain $\Omega_n$ is bounded, the maximum principle implies that $0\leq\psi_n\leq 1$ in $\Omega_n$.

Clearly, $\Omega$ is the limit domain of $\{\Omega_n\}_n$. 
Next we are going to prove that $\{u_n\}_n$ converges uniformly on compact sets of $\Omega$.
Let $D$ be a bounded subdomain of $\Omega$. For $n$ large enough, we have $D\subset\Omega_n$ so
we can restrict $u_n$ to $D$ and apply the results exposed in Section~\ref{secprelim}.

Firstly, we are going to prove there are no divergence lines for $\{u_n\}_n$.
Suppose by contradiction that there exists a divergence line $L$. Since $0\leq\psi_n\leq 1$, we deduce from Lemma~\ref{lem2} that $L$ must have
length no bigger than one. 
Thus taking a larger domain $D$ if necessary and using Lemma~\ref{lem1}, we obtain $L$ has to be a segment $[p_i,p_j]$.
If $i$ and $j$ have the same parity, then $\psi_n(p_i)=\psi_n(p_j)$ so $L$ has length zero (again using Lemma~\ref{lem2}),
which is absurd. If $i$ and $j$ have different parity, then $|\psi_n(p_i)-\psi_n(p_j)|=1$, so
$L$ is a chord of length one between $p_i$ and $p_j$. However, this is impossible on a non-special domain (see Proposition~\ref{prop-special}).

Hence there exists a subsequence of $\{u_n\}_n$ which converges on compact subsets of $D$.
Taking an exhaustion of $\Omega$ by bounded subdomains and using a diagonal process, we obtain a subsequence
of $\{u_n\}_n$ converging on compact subsets of $\Omega$ to a solution $u$ of the minimal graph equation.
By Lemma~\ref{lem4}, $u$ takes the marked values $\pm\infty$ on $\partial\Omega$. Its conjugate function
$\psi_{u}$ is the limit of $\{\psi_n\}_n$, hence $0\leq \psi_{u}\leq 1$ in $\Omega$.
By Lemma~\ref{lem5}, $u$ is the unique solution to the Jenkins-Serrin problem with bounded
conjugate function, and Proposition~\ref{prop1} follows. 
(Note we deduce from uniqueness that the whole sequence $\{u_n\}_n$ converges to $u$).
\end{proof}

\begin{remark}
In general, if $u$ is a solution to the Jenkins-Serrin problem on $\Omega$, then its conjugate function
$\psi_u$ satisfies $0\leq \psi_u\leq 1$ on the boundary of $\Omega$. However, if $\Omega$ is not contained
in a strip, $\psi_u$ might very well be unbounded, in which case we could not use the maximum principle
to guarantee that $0\leq \psi_u\leq 1$ in $\Omega$. This is why we took special care to construct $\Omega_n$ and $u_n$
in such a way that $\psi_n$ is bounded.
\end{remark}

\begin{figure}\begin{center}
\epsfysize=5cm \epsffile{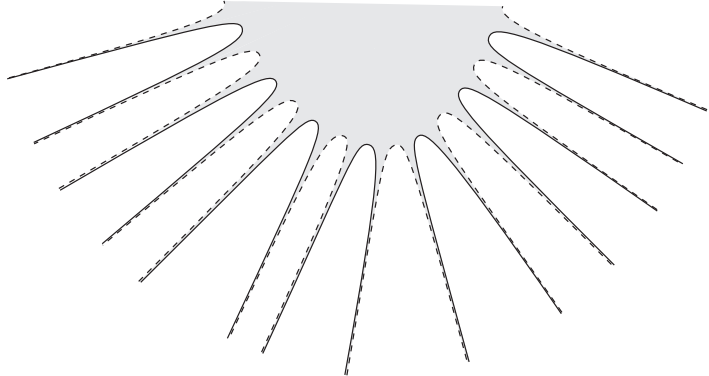}
\end{center}\caption{The domain on which the conjugate surface $M^*$ is a graph}
\end{figure}

Let $M$ be the graph of $u$ on $\Omega$. It is a minimal surface bounded by infinitely many vertical
straight lines above the vertices of $\Omega$.
Let $n_i$ be the normal to the open edge $(p_i,p_{i+1})$ pointing outwards $\Omega$.
Along the edge $(p_i,p_{i+1})$, the downward pointing normal to $M$ converges to $(-1)^i n_i$.
Let $M^*$ be the conjugate minimal surface of $M$. 
Since $0\leq \psi_u\leq 1$, $M^*$ is included in the slab $\{0\leq x_3\leq 1\}$.
Moreover, $\psi_u=0$ (resp. 1) at $p_i$ when $i$ is even (resp. odd), so the vertical line above $p_i$ on $M$
corresponds in the conjugate surface $M^*$ to an infinite horizontal symmetry curve lying on
the plane $\{x_3=0\}$ (resp. $\{x_3=1\}$).
The normal along this curve rotates from $(-1)^{i-1}n_{i-1}$ to $(-1)^i n_i$.
Finally, $M$ is a graph on a convex domain, thus $M^*$ is also a graph
on a (non convex) domain by a theorem of R. Krust.

Extending $M^*$ by symmetry with respect to the horizontal planes at integer heights, we obtain a
complete properly embedded singly periodic minimal surface ${\cal M}$ with period $(0,0,2)$.
It is easy to see that he quotient of ${\cal M}$ by its period has genus $0$, infinitely many ends and one limit end. 
To finish Theorem~\ref{mainth}, it only remains to prove that ${\cal M}$ has bounded curvature. 
Next section is devoted to that.

\begin{figure}
\epsfysize=6cm\begin{center}\epsffile{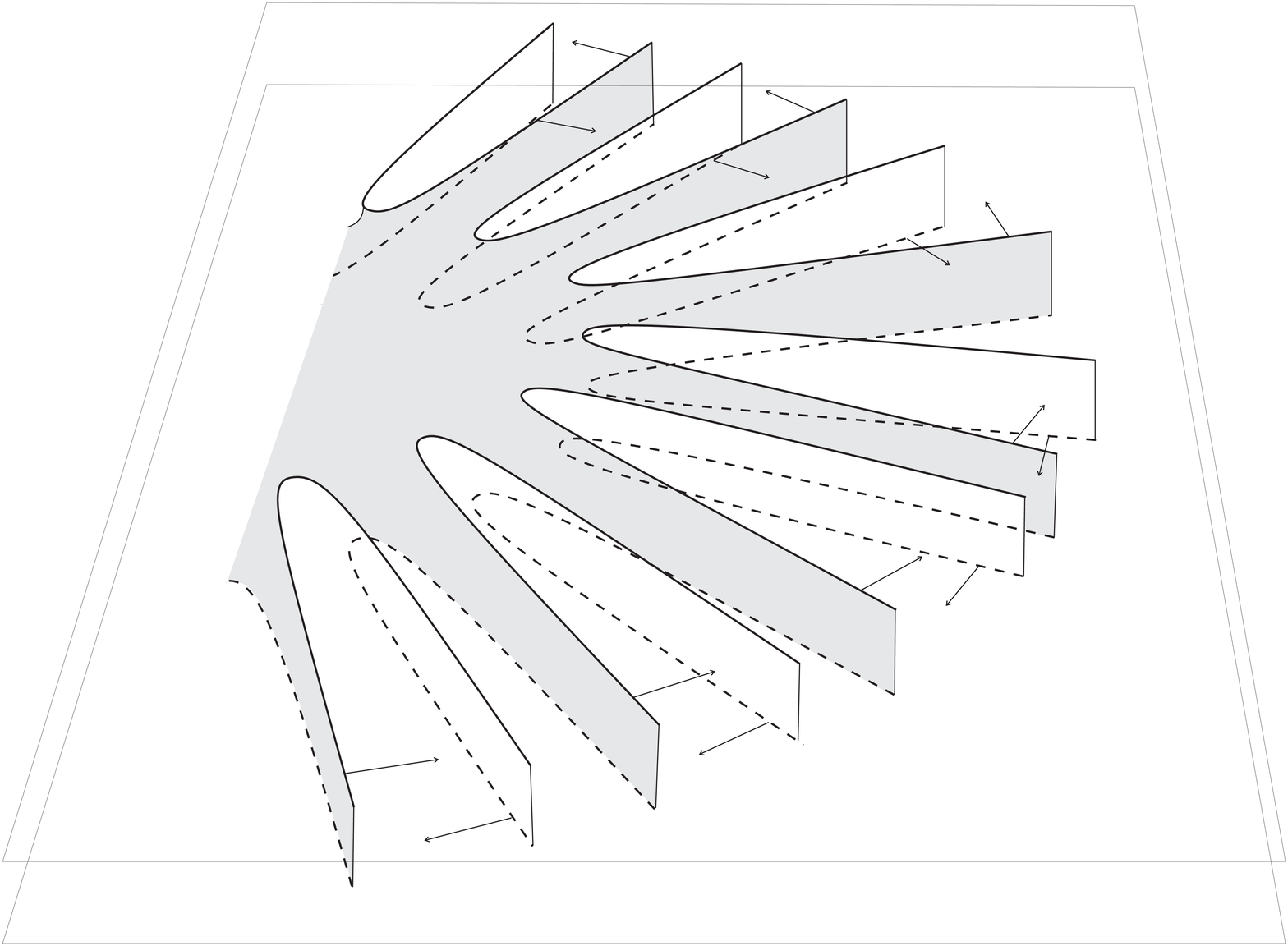}
\end{center}\caption{A sketch of the conjugate surface $M^*$} \end{figure}

\section{Bounded curvature}
\label{secBounded}
In the previous section, we constructed, for each unbounded convex polygonal domain $\Omega$ as
in the introduction, a minimal surface ${\cal M}$. 
Since $\Omega$ is convex and $\partial\Omega$ has unitary edges, 
there exists $r_0\in(0,1)$ such that the following holds: 
for each vertex $p$ of $\Omega$, the disk $D(p,r_0)$ intersects $\partial\Omega$ only along the two edges
with endpoint $p$; or equivalently, $D(p,r_0)\cap\Omega$ is a circular sector (whose angle is at most $\pi$, by convexity).

\begin{proposition}
\label{propBounded}
There exists a constant $C$ (independent of $\Omega$) such that the curvature $K$ of ${\cal M}$
is bounded by $C/r_0^2$.
\end{proposition}
\begin{proof}
As $M$ and $M^*$ (defined as in Section~\ref{secJSproblem}) are isometric, it suffices to bound the curvature of~$M$. Recall that
$M$ is the graph of a solution $u$ of (\ref{eqmin}) on $\Omega$ with boundary values $\pm\infty$ disposed alternately.
By Lemma \ref{lemmaVertex} below,
there exists $\delta>0$ such that, for each vertex $p$ of $\partial\Omega$, 
we have $W>10$ in $D(p,\delta r_0)\cap \Omega$, where $W=\sqrt{1+|\nabla u|^2}$. 

Given $s>0$, let $\Sigma_{p,s}$ be the graph of $u$ over $D(p,s)\cap\Omega$, completed by symmetry
about the vertical line above $p$. 
Note that $\Sigma_{p,s}$ is no more a graph.
The third coordinate of the normal to the graph of $u$ equals $-1/W$, thus the image of $\Sigma_{p,\delta r_0}$
by the Gauss map is contained in the domain $\S^2\cap\{|x_3|<1/10\}$. 
Since the area of this spherical domain is certainly less than $2\pi$, $\Sigma_{p,\delta r_0}$ is stable by a theorem of Barbosa and do Carmo~\cite{bc1}.
The distance from a point in $\Sigma_{p,\delta r_0/2}$ to the boundary of $\Sigma_{p,\delta r_0}$ is greater
than $\delta r_0/2$, hence a theorem by Schoen~\cite{sc3} assures that the absolute curvature of $\Sigma_{p,\delta r_0/2}$ is bounded above by $c/(\delta r_0/2)^2$, for some universal constant $c$.

On the other hand, $M$ is stable (because it is a graph) and bounded by the vertical lines above the vertices of $\partial\Omega$.
If we consider a point in $M\setminus \left(\bigcup_p \Sigma_{p,\delta r_0/2}\right)$, its distance to the boundary of $M$
is greater than $\delta r_0/2$, so by Schoen's theorem again, the absolute curvature at this point is bounded above 
by $c/(\delta r_0/2)^2$. This finishes the proof of Proposition~\ref{propBounded}.  
\end{proof}

\medskip
We now prove a useful gradient estimate from below in a neighborhood of the vertices of circular sectors, 
used in the previous proof.
Given $0<\alpha<2\pi$, let $U_{\alpha}$ be the circular sector domain of radius one and angle $\alpha$, 
defined in polar coordinates by $\{0<r<1,\ 0<\theta<\alpha\}$.

\begin{lemma}
\label{lemmaVertex}
Given $C>0$, there exists $\delta>0$ (independent of $\alpha$) such that the following is true:
Let $u$ be a solution of (\ref{eqmin}) on $U_{\alpha}$ such that $u=+\infty$ on the segment $\{\theta=0\}$ 
and $u=-\infty$ on $\{\theta=\alpha\}$.
Let $\psi_u$ be the conjugate function of $u$ normalized so that $\psi_u(0)=0$, and assume $\psi_u\geq 0$ in $U_{\alpha}$.
Then, $|\nabla u|\geq C$ in $U_{\alpha}\cap D(0,\delta)$.
\end{lemma}
\begin{proof}
Assume by contradiction that the lemma is not true.
Then there exist sequences $\{\alpha_n\}_n$, $\{u_n\}_n$ and $\{x_n\}_n$ such that:
\begin{itemize}
\item $u_n$ is a solution of (\ref{eqmin}) on $U_{\alpha_n}$ with boundary values $u_n=+\infty$ on $\{\theta=0\}$ 
and $u_n=-\infty$ on $\{\theta=\alpha_n\}$;
\item $x_n\in U_{\alpha_n}$, $x_n\to 0$;
\item $|\nabla u_n(x_n)|<C$ for every $n$.
\end{itemize}
We deduce from Lemma~1 of Jenkins and Serrin~\cite{jes1} that there exists ${\varepsilon>0}$ (depending on $C$) 
such that $|\nabla u_n|\geq C$ in the sector defined in polar coordinates by $\{0<r<1/8,\ 0<\theta<\varepsilon\}$, 
and the same holds in $\{0<r<1/8,\ \alpha_n-\varepsilon<\theta<\alpha_n\}$.
The existence of $x_n$ implies that $\alpha_n\geq 2\varepsilon$.
After passing to a subsequence, we can assume that $\{\alpha_n\}_n$ converges to some 
$\alpha_{\infty}\in [2\varepsilon,2\pi]$.

Let $\lambda_n=1/|x_n|$, $\widetilde{U}_n=\lambda_n U_{\alpha_n}$ and $\widetilde{x}_n=\lambda_n x_n$.
Define $\widetilde{u}_n(x)=\lambda_n u_n(x/\lambda_n)$, for every $x\in\widetilde{U}_n$  
(the graph of $\widetilde{u}_n$ coincides with the graph of $u_n$ scaled by $\lambda_n$). 
Passing to a subsequence, we have $\widetilde{x}_n\to\widetilde{x}_{\infty}\in\{r=1,\theta=\theta_{\infty}\}$, 
with $\varepsilon\leq\theta_{\infty}\leq\alpha_{\infty}-\varepsilon$. 
By slightly rotating the domains $\widetilde{U}_n$, we may assume that $\widetilde{x}_n=\widetilde{x}_{\infty}$.
The limit domain of the (rotated) $\widetilde{U}_n$ is the unbounded sector 
$\widetilde{U}_{\infty}=\{0<r<\infty,0<\theta<\alpha_{\infty}\}$.

Since $|\nabla \widetilde{u}_n(\widetilde{x}_n)|=|\nabla u_n(x_n)|\leq C$,
the convergence domain ${\cal B}(\widetilde u_n)$ of the sequence $\{\widetilde{u}_n\}_n$ contains
the point $\widetilde{x}_\infty$.
Let $D$ be the component of ${\cal B}(\widetilde u_n)$ which contains $\widetilde{x}_\infty$, 
and let $\widetilde u_\infty$ be the limit of a subsequence of $\{\widetilde{u}_n-\widetilde{u}_n(\widetilde{x}_{\infty})\}_n$ on $D$. 
The boundary of $D$ consists of divergence lines of $\{\widetilde{u}_n\}_n$, 
so let us see what divergence lines are possible.
Using Lemma \ref{lem2} and the fact that the conjugate function $\widetilde{\psi}_n$ of $\widetilde{u}_n$ is positive, 
we obtain that a divergence line cannot be a complete line in the plane. 
And by Lemma~\ref{laurent}-{\it(i)}, we know that a divergence line can only meet the boundary of $\widetilde{U}_{\infty}$ at the origin.
Hence the only possible divergence lines for $\{\widetilde{u}_n\}_n$ are half lines with endpoint at the origin, 
and $D$ is an unbounded sector bounded by two half lines $L_1=\{0<r<\infty,\theta=\beta_1\}$ and  
$L_2=\{0<r<\infty,\theta=\beta_2\}$, with $0\leq \beta_1<\beta_2\leq\alpha_{\infty}$. 

If $\beta_1>0$, we conclude from lemmas \ref{lem2} and \ref{lem4} that $\widetilde{u}_n=+\infty$ on $L_1$. 
If $\beta_1=0$ then the same is true using Lemma~\ref{laurent}-{\it (ii)}. 
In the same way, we obtain that $\widetilde{u}_n=-\infty$ on $L_2$.
It is proven in \cite{mazet4}, Proposition~2 (or Proposition~4, when $D$ is a sector of angle $2\pi$), 
that the Jenkins-Serrin problem on this unbounded sector has no solution.  This contradiction proves the lemma.
\end{proof}

\section{Asymptotic behavior}
\label{secAsymptotic}
Assume that $\Omega$ is {\em not} contained in a strip. We will prove that the surface ${\cal M}$ we constructed in the section \ref{secJSproblem}
is asymptotic, in a sense that we will explain,
to two doubly periodic Scherk minimal surfaces.
When $\Omega$ is contained in a strip, it may be proven that the surface is asymptotic
to a Toroidal Halfplane Layer~\cite{ka4,mr3,mrod1}, which is a doubly periodic minimal surface with parallel Scherk-type ends 
(they have been classified in \cite{PeRoTra1} as the only properly embedded doubly periodic minimal surfaces with genus $1$ and a finite number of ends in the quotient).
The argument is similar, although a little more involved. Thus we will only consider here the
case where $\Omega$ is not contained in a strip.

Let $a_n=p_{n+1}-p_n\in\S^1$. Since $\Omega$ is convex, the limits $a_{\infty}=\lim_{n\to\infty}a_n$ and
$a_{-\infty}=\lim_{n\to -\infty}a_n$ exist.
Let $a_{\infty}^{\perp}$ be the unitary horizontal vector orthogonal
to $a_{\infty}$, chosen so that
$q_n=p_n+a_{\infty}^{\perp}\in\Omega$ for $n$ large enough.
Let $Q_n=(q_n,u(q_n))$ be the corresponding point in $M$ and
$Q_n^*$ be the corresponding point in $M^*$.
\begin{proposition}
\label{propAsymptotic}
When $n\to\infty$, ${\cal M}$ translated by $-Q_{2n}^*$ converges to a doubly periodic Scherk minimal surface with
periods $(0,0,2)$ and $(2 a_{\infty}^{\perp},0)$.
The convergence is smooth convergence on compact subsets of $\R^3$.
A similar statement holds when $n\to -\infty$.
\end{proposition}
\begin{proof}
Let $\widetilde{\Omega}_n$ be the domain $\Omega$ translated by $-p_{2n}$.
The limit domain of the sequence $\{\widetilde{\Omega}_n\}_n$
is a half plane $\widetilde{\Omega}_{\infty}$ bounded by the line $\mbox{Span}(a_{\infty})$.
Hence it is natural to study the corresponding Jenkins-Serrin problem on this half plane.

Without loss of generality, we may assume that $\widetilde{\Omega}_{\infty}$ is the half plane
$\{x_2\geq 0\}$, so $a_{\infty}^{\perp}=(0,1)$, and the boundary data is
 $+\infty$ on $[\widetilde{p}_i,\widetilde{p}_{i+1}]$ if $i$ is even, 
and $-\infty$ if $i$ is odd, where $\widetilde{p}_i=(i,0)$.
This Jenkins-Serrin problem on $\widetilde{\Omega}_{\infty}$ has the following explicit solution:
Let $U$ be the half band $\{0\leq x_1\leq 1,\ x_2\geq 0\}$, with boundary data $+\infty$ on
the horizontal segment and $0$ on the vertical half lines.
A piece of singly periodic Scherk minimal surface, rotated so that its period is $(2,0,0)$, 
solves the Jenkins-Serrin problem on $U$. 
Extending by symmetry, we obtain a solution to the Jenkins-Serrin problem on $\widetilde{\Omega}_{\infty}$. 
Let us call $u_S$ this solution and $S$ its graph.
The conjugate minimal surface $S^*$ of $S$ is a piece of doubly periodic Scherk minimal surface, 
rotated so that its periods are $(0,2,0)$ and $(0,0,2)$, such piece lying in the slab $\{0\leq x_3\leq 1\}$. 
In particular, the conjugate function $\psi_S$ of $u_S$ is bounded. 
By Lemma~\ref{lem5}, $u_S$ is the unique solution to the Jenkins-Serrin problem on
$\widetilde{\Omega}_{\infty}$ with bounded conjugate function.

\begin{claim}
\label{claim1}
Let $\widetilde{u}_n(x)=u(x+p_{2n})$, for all $x\in\widetilde{\Omega}_n$. 
Then, $\{\widetilde{u}_n-\widetilde{u}_n(0,1)\}_n$ converges on compact subsets of
$\widetilde{\Omega}_{\infty}$ to $u_S$.
\end{claim}

Since $\widetilde{\Omega}_{\infty}$ is a half-plane, 
every divergence line $L$ of $\{\widetilde{u}_n\}_n$ can be extended as far as we want in at least one direction.
However, since $0\leq\widetilde{\psi}_n\leq 1$, Lemma~\ref{lem2} says that $L$ has length at most one. 
Hence we deduce there are no divergence lines, and the sequence $\{\widetilde{u}_n\}_n$ converges on compact
subsets of $\widetilde{\Omega}_{\infty}$ to a solution $\widetilde{u}_{\infty}$ of (\ref{eqmin}).
By Lemma~\ref{laurent}-{\it (ii)}, this solution has boundary value $+\infty$
on the open segment $(\widetilde{p}_i,\widetilde{p}_{i+1})$ if $i$ is even, 
and $-\infty$ if $i$ is odd, so it solves the Jenkins-Serrin problem on $\widetilde{\Omega}_{\infty}$.
By uniqueness (Lemma \ref{lem5}), we have $\widetilde{u}_{\infty}=u_S$, and Claim~\ref{claim1} is proven.

\medskip
Let us return to the proof of Proposition~\ref{propAsymptotic}.
It is tempting to say that $M-Q_{2n}$ converges to a singly periodic Scherk minimal surface so
$M^* -Q_{2n}^*$ converges to a doubly periodic Scherk minimal surface.
There are several problems with this approach. First we would need a notion
of convergence for surfaces with boundary. 
The main problem is that the convergence of $\widetilde{u}_n$ to $u_S$, Claim~\ref{claim1},
only holds on compact subsets of $\widetilde{\Omega}_\infty$, so it does not say us what happens in a neighborhood of
the vertical lines.
For this reason, we argue as follows:
${\cal M}$ is a properly embedded minimal surface in $\R^3$ whose absolute curvature is bounded by
some constant $c$ (Proposition~\ref{propBounded}).
By the Regular Neighborhood Theorem, or ``Rolling Lemma'' 
(firstly proven by A.~Ros~\cite{ros5}, Lemma 4, for properly embedded minimal surfaces in $\R^3$ with  finite total curvature, 
and generalized to properly embedded minimal surfaces with bounded curvature by 
Meeks and Rosenberg~\cite{mr7}, Theorem 5.3), 
${\cal M}$ has an embedded tubular neighborhood of radius $1/\sqrt{c}$. 
In particular, we have local area bounds, namely, the area of ${\cal M}$ 
inside balls of radius $1/\sqrt{c}$ is bounded by some constant. 
By standard result, a subsequence of $\{{\cal M} - Q_{2n}^*\}_n$
converges with finite multiplicity, on compact subsets of $\R^3$, to a limit minimal
surface. From Claim~\ref{claim1}, the limit must be a Scherk doubly periodic surface, and the
multiplicity is one. So the whole sequence converges to a Scherk doubly periodic surface.
This concludes the proof of Proposition~\ref{propAsymptotic}. 
\end{proof}

\section{Appendix}
For completeness, we prove in this section that, among all the bounded convex unitary polygonal
domains, the ones that fail to satisfy the hypothesis of the theorem of Jenkins and Serrin
are the special domains.

Let $\Omega$ be a bounded convex polygonal domain, with sides marked alternately by 
$+\infty$ and $-\infty$, and ${\cal P}$ be any polygonal subdomain of $\Omega$ (this means that its
vertices are vertices of $\Omega$). 
Denote by $\alpha$ (resp. $\beta$) the total length of the edges of ${\cal P}$ which are edges of $\Omega$ with
mark $+\infty$ (resp. $-\infty$), and call $\gamma$ the perimeter of ${\cal P}$. 
The domain $\Omega$ satisfies the hypothesis of the theorem of Jenkins and Serrin 
(and so one can solve the Jenkins-Serrin problem on $\Omega$) if and only if 
$2\alpha<\gamma$ and $2\beta<\gamma$ for each strict subpolygon ${\cal P}$ of $\Omega$, 
and $\alpha=\beta$ when ${\cal P}=\Omega$.

Consider a convex polygonal domain $\Omega$ as above, and suppose all its edges have length one. 
Label its vertices $p_1,\cdots,p_{2n}$ so that $[p_1,p_2]$ is marked with $-\infty$ 
(so $[p_i,p_{i+1}]$ is marked with $+\infty$ if $i$ is even and with $-\infty$ if $i$ is odd, with the convention $p_{2n+1}=p_1$).
We say $p_i$ is an even vertex if $i$ is
even and an odd vertex if $i$ is odd.
We will refer as a {\it chord} to a straight segment 
that joints two different non consecutive vertices of $\Omega$.
\begin{proposition}
\label{prop-special}
The following statements are equivalent:
\begin{itemize}
\item[(i)] $\Omega$ is not a special domain.
\item[(ii)] Every chord from an even vertex to an odd vertex has
length greater than 1.
\item[(iii)] $\Omega$ satisfies the hypotheses of the theorem of
Jenkins and Serrin.
\end{itemize}
\end{proposition}

Before proving Proposition~\ref{prop-special}, let us recall the following elementary result proven in~\cite{PeTra1}, Lemma 5.2.

\begin{lemma}\label{lem6}
Let $ABCD$ be a convex quadrilateral such that $|BC|=|AD|$ and 
$\widehat{A}+\widehat{B}\leq\pi$, where $\widehat{A}$ means the interior angle at $A$.
Then $|CD|\leq |AB|$, with equality if and only if $ABCD$ is a parallelogram.
\end{lemma}

\begin{proof}
Let us see that $(i)\Rightarrow (ii)$. Arguing by contraposition, we must prove that
if $\Omega$ has a chord of length $\leq 1$ from an even vertex to an odd vertex, then
$\Omega$ is special.
Let $C$ be such a chord. It divides $\Omega$ into two convex domains, $\Omega_1$ and $\Omega_2$.
For one of them, let us say $\Omega_1$,
the sum of the inner angles at the endpoints of $C$ is $\leq \pi$. We may rename the vertices
of $\Omega$, without changing their parity, so that $C$ is the segment $[p_1,p_{2r}]$ and
the vertices on the boundary of $\Omega_1$ are $p_1\cdots,p_{2r}$.
Lemma~\ref{lem6} assures $|p_2p_{2r-1}|\leq |p_1p_{2r}|\leq 1$; 
and by induction, we obtain $|p_ip_{2r+1-i}|\leq 1$ for all $1\leq i\leq r$ 
(here we use that the sum of the inner angles remains $\leq\pi$ because $\Omega$
is convex).
But $[p_r,p_{r+1}]$ is an edge on the boundary of $\Omega$, so $|p_rp_{r+1}|=1$.
Hence equality holds everywhere, and all quadrilaterals $p_ip_{i+1}p_{2r-i}p_{2r+1-i}$
are parallelograms. 
Since $\Omega_1$ is convex, $p_1p_r p_{r+1}p_{2r}$ is a parallelogram (with
two sides of length $1$ and two sides of length $r-1$).
Hence the sum of the inner angles of $\Omega_1$ at the endpoints of $C$ is $\pi$, so the
sum of the inner angles of $\Omega_2$ at the same points is $\leq \pi$. Applying the same
argument to $\Omega_2$, we obtain that $\Omega_2$ is also a parallelogram, so $\Omega$
is a special domain.

Let us see that $(ii)\Rightarrow (iii)$.
Let ${\cal P}$ be a strict subpolygon of $\Omega$. Let us orient 
$\partial \Omega$ and $\partial {\cal P}$ as boundaries of
$\Omega$ and ${\cal P}$. Note that for an edge in $\partial {\cal P}
\cap \partial \Omega$, both orientations are the same. Let us prove that
$2\alpha<\gamma$. This is clearly true if $\partial {\cal P}$ contains no edge
marked $+\infty$. Let $[p_{2i},p_{2i+1}]$ be an edge on the boundary of ${\cal P}$
marked $+\infty$. Let $[p_{2j},p_{2j+1}]$ be the next edge on the boundary of ${\cal P}$
marked $+\infty$, when traveling along the boundary in the direction given by its orientation.
Let $C$ be the part of $\partial{\cal P}$ between $p_{2i+1}$ and $p_{2j}$.
If $C$ contains an edge marked $-\infty$ then $|C|\geq 1$. Else $C$ contains only chords. Since $C$ connects an odd vertex
with an even vertex, at least one of its chords goes from an odd vertex to an even vertex,
so $|C|>1$. Hence the part of the boundary of ${\cal P}$
between two edges marked $+\infty$ always has
length $\geq 1$, with strict inequality for at least one of them (else ${\cal P}=\Omega$).
Hence $2\alpha<\gamma$.
The proof of $2\beta<\gamma$ is exactly the same, exchanging the roles of $+\infty$ and
$-\infty$.

Finally, $(iii)\Rightarrow (i)$ is obvious: a special domain $\Omega$ does not satisfy
the hypothesis of Jenkins and Serrin, because if ${\cal P}$ is a rhombus then $2\alpha=\gamma$ (or
$2\beta=\gamma$).  
\end{proof}

\bibliographystyle{plain}

\end{document}